\documentclass[12pt,noamsfonts]{amsart}

\usepackage{amsfonts}
\usepackage{amssymb}
\usepackage{amscd}

\newcommand{\marginlabel}[1]%
  {\mbox{}\marginpar{\raggedleft\hspace{0pt}\bfseries\sf#1}}

\def\NN{{\mathbb N}}

\def\RR{{\mathbb R}}
\def\QQ{{\mathbb Q}}

\def\cI{\mathcal{I}}

\def\cO{\mathcal{O}}

\def\bm{\mathbf{m}}
\def\ba{\mathbf{a}}
\def\bb{\mathbf{b}}
\def\bc{\mathbf{c}}

\DeclareMathOperator{\Spec}{Spec}

\DeclareMathOperator{\vol}{vol}
\DeclareMathOperator{\lc}{lc}
\DeclareMathOperator{\lin}{in}
\DeclareMathOperator{\ee}{e}
\DeclareMathOperator{\length}{l}

\newtheorem{lemma}{Lemma}[section]
\newtheorem{theorem}[lemma]{Theorem}
\newtheorem{corollary}[lemma]{Corollary}

\theoremstyle{definition}
\newtheorem{definition}[lemma]{Definition}
\newtheorem{remark}[lemma]{Remark}

\newtheorem{example}[lemma]{Example}
\newtheorem{question}[lemma]{Question}

\theoremstyle{remark}
\newtheorem*{remark*}{Remark}
\newtheorem*{note*}{Note}

\begin{document}

\title{On multiplicities of graded sequences of ideals}

\author[M. Musta\c{t}\v{a}]{Mircea~Musta\c{t}\v{a}}
\address{Isaac Newton Institute for Mathematical Sciences,
20 Clarkson Road, Cambridge CB3 0EH, England}
\email{{\tt mirceamustata@yahoo.com}}

\thanks{This research was conducted during the period the author
served as a Clay Mathematics Institute Long-Term Prize Fellow}

\begin{abstract}
We generalize a result from \cite{els}, proving that for an arbitrary
graded sequence of zero-dimensional ideals, the multiplicity of the sequence
is equal to its volume. This is done using a deformation to monomial ideals.
As a consequence of our result, we obtain a formula which computes the
multiplicity of an ideal $I$ in terms of the multiplicities of the initial
monomial ideals of the powers $I^m$. We use this to give a new proof 
of the inequality between multiplicity and the log canonical threshold
from \cite{dfelm}.
\end{abstract}

\maketitle
\markboth{M. MUSTA\c T\v A}
{MULTIPLICITIES OF GRADED SEQUENCES OF IDEALS}

\section*{Introduction}

Let $R$ be a regular local ring of dimension $n$. A graded sequence
of ideals in $R$ is a set of ideals $\ba_{\bullet}=
\{\ba_m\}_{m\in\NN}$ such that for all $p$, $q$, we have
$\ba_p\cdot
\ba_q\subseteq\ba_{p+q}$. The trivial example is given by the
powers of a fixed ideal. A more interesting case is that of the
symbolic powers of a given ideal.
Geometric examples arise as follows:
let $X$ be a smooth variety and $L$ a line bundle on $X$. If
$R=\cO_{X,Z}$, for some irreducible, closed $Z\subseteq X$,
and if $\ba_m$ defines in $R$ the base locus of the complete linear system
$|L^m|$, then $\ba_{\bullet}$ is a graded sequence of ideals.

In \cite{els}, Ein, Lazarsfeld and Smith have introduced the volume
and the multiplicity of a graded sequence $\ba_{\bullet}$
of zero-dimensional ideals. The volume is
defined by 
$$\vol(\ba_{\bullet}):=\limsup_{m\to\infty}\frac{n!\cdot \length(R/\ba_m)}{m^n},$$
while the multiplicity is given by
$$\ee(\ba_{\bullet}):=\lim_{m\to\infty}\frac{\ee(\ba_m)}{m^n},$$
where for a zero-dimensional ideal $I$, we denote by $\ee(I)$ the
Hilbert-Samuel multiplicity of $R$ along $I$. It was proved in
\cite{els} that under a certain condition on $\ba_{\bullet}$
(see below for details), 
we have $\ee(\ba_{\bullet})=\vol(\ba_{\bullet})$. 
The proof was based on the theory of asymptotic multiplier ideals,
so it required the restriction to characteristic zero.
This result was applied to study Abhyankar valuations.

Our main result is that the equality $\ee(\ba_{\bullet})=\vol(\ba_{\bullet})$
holds for an arbitrary graded sequence of zero-dimensional ideals. The main
idea of the proof is to reduce the assertion to the case of a graded sequence
of monomial ideals. In particular, the proof is purely algebraic, so it 
applies to any regular ring containing a field (of arbitrary characteristic).

As a byproduct of our proof, we obtain a useful result, which is new
even in the case of one ideal.
Suppose that $I\subset R=K[X_1,\ldots,X_n]$ is an ideal supported at the origin.
Fix a monomial order and let $\ba_m=\lin_>(I^m)$. We obtain as a corollary
the following formula:
\begin{equation}\label{formula}
\ee(I)=\lim_{m\to\infty}\frac{\ee(\ba_m)}{m^n}.
\end{equation}

We apply this result to give an easier proof of the inequality 
from \cite{dfelm} between
the log canonical threshold $\lc(I)$ and the multiplicity
$\ee(I)$ (see \cite{dfelm} for motivation
in the context of birational geometry). More precisely, we show that
if $I$ is a zero-dimensional ideal in an $n$-dimensional
local ring of a smooth complex variety,
then we have
$$\ee(I)\geq \frac{n^n}{\lc(I)^n}.$$
The main point is that $(\ref{formula})$ reduces the statement to the
case of a monomial ideal, when the inequality is easy to check.

\bigskip

We explain now in more detail the idea of the proof and how it relates 
to the approach in \cite{els} based on multiplier ideals. A few words about
these ideals: they  have been introduced in the analytic setting by Demailly,
Nadel and Siu, but they have found striking applications 
in algebraic geometry, as well, in the work of Ein and
Lazarsfeld and Kawamata (see, for example, \cite{el}, \cite{kawamata}
and \cite{siu}). Multiplier ideals, and especially their asymptotic version,
have turned out to be a powerfull tool also in the study 
of graded sequences of ideals in an algebraic setting
(see \cite{els2}, where they are used to relate the symbolic powers
and the usual powers of an ideal).

The asymptotic ideals of a graded sequence $\ba_{\bullet}$ form a family 
$\bb_{\bullet}=\{\bb_m\}_{m\in\NN^*}$ of ideals 
such that $\ba_m\subseteq\bb_m$ for all $m$ and which satisfy a property
which is ``opposite''
 to the defining property of $\ba_{\bullet}$: for every $p$ and $q$,
$\bb_{p+q}\subseteq\bb_p\cdot\bb_q$ (this is the Subadditivity
Theorem of \cite{del}). 
In particular, for every $m$ and $p$, we have the inclusions
$$\ba_m^p\subseteq\ba_{mp}\subseteq\bb_{mp}\subseteq\bb_m^p.$$
Very loosely, one could say that $\bb_{\bullet}$ can be used to 
compensate the failure of $\ba_{\bullet}$ to be the sequence of 
powers of an ideal.

Suppose now that $\ba_{\bullet}$ is a graded sequence of zero-dimensional
ideals. One can define as above the invariants $\ee(\bb_{\bullet})$
and $\vol(\bb_{\bullet})$ and it is easy to see that
$\ee(\bb_{\bullet})\leq\vol(\bb_{\bullet})\leq\vol(\ba_{\bullet})\leq
\ee(\ba_{\bullet})$. The result in \cite{els} is that if the graded
sequence of colon ideals $\{\ba_m:\bb_m\}_m$ has multiplicity zero, then 
$\ee(\ba_{\bullet})=\ee(\bb_{\bullet})$. 
This condition is satisfied, for example,
by the sequence defined by an Abhyankar valuations, as in \cite{els},
or by the sequence defining the base loci of the powers of a big line bundle.

We show that if $\ba_{\bullet}$ consists of monomial ideals
in a polynomial ring over a field, then
$\ee(\ba_{\bullet})=\ee(\bb_{\bullet})$. 
Note that
in this case, multiplier ideals can be 
introduced in terms of Newton polyhedra 
(this is a result from \cite{howald})
and the main property, subadditivity, can be easily proved directly.
What we show, in fact,
is that a sequence of ideals closely related to $\ba_{\bullet}$
satisfies the criterion in \cite{els} and that this is enough to give
$\ee(\ba_{\bullet})=\ee(\bb_{\bullet})$.

By deforming an arbitrary graded sequence of ideals to a monomial sequence,
we deduce that $\ee(\ba_{\bullet})=\vol(\ba_{\bullet})$ in general. Moreover,
we show that $\ee(\ba_{\bullet})=0$ if and only if there is 
$p\in\NN^*$ such that $\ba_q\subseteq {\bf m}^{[q/p]}$, for all $q$,
where ${\bf m}$
is the maximal ideal of $R$. When we 
are in a situation where multiplier ideals are defined, this means
that $\ee(\ba_{\bullet})=0$ if and only if $\ee(\bb_{\bullet})=0$.

We do not know whether we 
always have 
$\ee(\ba_{\bullet})=\ee(\bb_{\bullet})$
(assuming, of course, that $\bb_{\bullet}$ is defined).
 We show, however, that if instead of 
multiplicity we consider the log canonical threshold, then we have equality:
$\lc(\ba_{\bullet})=\lc(\bb_{\bullet})$ (see \S 3 for definitions).

\smallskip

A few words about the structure of the paper: in the first section we
discuss the definition of multiplicity and volume, and reduce the statement
of the main theorem to the case of monomial ideals. In 
the next section, we treat
monomial ideals: we discuss 
asymptotic multiplier ideals,
and prove that $e(\ba_{\bullet})=e(\bb_{\bullet})$ in this case. 
The last section
applies the previous ideas to discuss another invariant, the log canonical threshold
for graded sequences of ideals. In particular, we prove that $\lc(\ba_{\bullet})=
\lc(\bb_{\bullet})$. We also apply our main result to deduce the inequality
involving the multiplicity and the log canonical threshold.

\bigskip

\subsection{Acknowledgements}
This work was motivated by the article \cite{els}. We are grateful
to Lawrence Ein, Robert Lazarsfeld and Karen Smith for providing us an early
version of their manuscript, without which the present paper would not exist.
Moreover, discussions with Lawrence Ein on this material have been
of invaluable help.

Work on the paper has been done while the author was visiting
Universit\'{e} de Nice-Sophia Antipolis and Isaac Newton Institute
for Mathematical Sciences. We are grateful to both institutions for
providing an excellent working environment.

\section{Volume versus multiplicity}

Let $(R,\bm)$ be a regular local ring containing a field, with
$\dim(R)=n$. 
Recall that a graded sequence of ideals in $R$ is a set of ideals
$\ba_{\bullet}=\{\ba_m\}_{m\in\NN}$ such that $\ba_p\cdot\ba_q
\subseteq\ba_{p+q}$ for every $p$, $q\in\NN$. 

The following definition for the volume of $\ba_{\bullet}$
appears in \cite{els}.

\begin{definition}\label{def_volume}
Suppose that $\ba_{\bullet}$ is a graded sequence of 
zero-dimensional ideals in $R$, i.e., $\dim\,(R/\ba_m)\leq 0$
for all $m\in\NN$.
 The volume of $\ba_{\bullet}$ is defined by
$$\vol(\ba_{\bullet}):=
\limsup_{m\to\infty}\frac{n!\cdot \length(R/\ba_m)}{m^n}.$$
\end{definition}

\begin{remark}
In \cite{els} one considers more general families of ideals, indexed by
an ordered semigroup $\Gamma$. However, one makes the 
additional assumption that
$\ba_p\subseteq \ba_q$ if $p\geq q$ in $\Gamma$. In this case, one can always
reduce the computation of volumes and multiplicities to families indexed
by $\NN$ (see \cite{els} for details).
\end{remark}

Our main result expresses the volume as a limit of Hilbert-Samuel
multiplicities of the ideals $\ba_m$. If $I$ is a zero-dimensional
ideal of $R$, we will denote by $\ee(I)$ the Hilbert-Samuel multiplicity
of $R$ along $I$. 
We start with the following 
easy lemmas which allow us to define the multiplicity of $\ba_{\bullet}$.

\begin{lemma}\label{ineq_mult}
For every graded sequence of zero-dimensional ideals $\ba_{\bullet}$,
and every $p$, $q\in\NN$, we have
$$\ee(\ba_{p+q})^{1/n}\leq \ee(\ba_p)^{1/n}+\ee(\ba_q)^{1/n}.$$
\end{lemma}

\begin{proof}
Since $\ba_{\bullet}$ is a graded system of ideals, we have
$\ba_p\cdot\ba_q\subseteq \ba_{p+q}$, hence $\ee(\ba_{p+q})
\leq \ee(\ba_p\cdot\ba_q)$. The assertion in the lemma follows from
this and Teissier's inequality (see \cite{teissier}):
$$\ee(\ba_p\cdot\ba_q)^{1/n}\leq \ee(\ba_p)^{1/n}+\ee(\ba_q)^{1/n}.$$
\end{proof}

\begin{lemma}\label{limit}
Let $\{\alpha_m\}_{m\in\NN}$ be a sequence of real numbers,
 with $\alpha_m\geq 0$
for every $m$.
If $\alpha_{p+q}\leq\alpha_p+\alpha_q$ for all 
$p$, $q\in\NN$, then $\lim_{m\to\infty}\frac{\alpha_m}{m}$
exists and it is equal to $\inf_{m\in\NN^*}\frac{\alpha_m}{m}$. 
\end{lemma}

\begin{proof}
Let $L=\inf_{m\in\NN^*}\frac{\alpha_m}{m}$ and $\epsilon>0$.
We show that $\frac{\alpha_q}{q}\leq L+\epsilon$ for
$q\gg 0$. 

We can find $m\geq 1$ such that $\alpha_m/m<L+\epsilon/2$.
For every integer $p$ with $0\leq p<m$ and every $k\in\NN$,
we have $\alpha_{km+p}\leq
k\alpha_m+\alpha_p$, hence
$$\frac{\alpha_{km+p}}{km+p}\leq\frac{km(L+\epsilon/2)+\alpha_p}
{km+p}.$$
When $k$ goes to infinity, the right hand side of the above inequality
goes to $L+\epsilon/2$. Hence for $k\gg 0$ we get 
$\frac{\alpha_{km+p}}{km+p}\leq L+\epsilon$. Since this holds
for every integer $p$, with $0\leq p<m$, we are done. 
\end{proof}

By combining the previous lemmas, we deduce the following

\begin{corollary}
If $\ba_{\bullet}$ is a graded system of zero-dimensional ideals,
then $\lim_{m\to\infty}\frac{\ee(\ba_m)}{m^n}$ exists and it is equal
to $\inf_{m\in\NN^*}\frac{\ee(\ba_m)}{m^n}$. 
\end{corollary}

\begin{definition}
If $\ba_{\bullet}$ is a graded sequence of zero-dimensional ideals,
then the multiplicity of $\ba_{\bullet}$ is defined by
 $$\ee(\ba_{\bullet}):=\lim_{m\to\infty}\frac{\ee(\ba_m)}{m^n}.$$
It is clear that $\ee(\ba_{\bullet})\in\RR_+$.
\end{definition}

For a real number $x$,
we denote by $[x]$ the integral part of $x$, i.e., the
largest integer $m$ such that $m\leq x$. The following is our main result.

\begin{theorem}\label{main}
Let $\ba_{\bullet}$ be a graded sequence of zero-dimensional ideals.
\begin{enumerate}
\item We have $\vol(\ba_{\bullet})=\ee(\ba_{\bullet})$.
\item $\vol(\ba_{\bullet})>0$ if and only if there is $q\in\NN^*$
such that $\ba_p\subseteq\bm^{[p/q]}$ for all $p\in\NN$.
\end{enumerate}
\end{theorem}

\begin{remark}
Over a field of characteristic zero, this was proved in \cite{els}
under the assumption that if $\bb_{\bullet}$ is the corresponding
sequence of asymptotic multiplier ideals, then $\ba_{\bullet}$
is close to $\bb_{\bullet}$ in a suitable sense (we recall the precise
statement in Lemma~\ref{key} below).
Under this extra hypothesis it is shown that, in fact, the multiplicity of
$\ba_{\bullet}$ can be computed as the multiplicity $\ee(\bb_{\bullet})$
of $\bb_{\bullet}$.
We do not know whether this also holds for
an arbitrary sequence $\ba_{\bullet}$, but we will show
in $\S 2$ that the assertion is true for graded sequences of monomial ideals.
 Note also that the 
second assertion in the above theorem can be interpreted as saying that
for arbitrary $\ba_{\bullet}$, we have $\ee(\ba_{\bullet})=
\ee(\bb_{\bullet})$
if one of these invariants is
zero (see \S 2 for details).
\end{remark}

The equality in Theorem~\ref{main} $(1)$ allows us to deduce a
generalization to volumes of Teissier's inequality for multiplicities.
 The idea is the same as in \cite{els}, but now
we get the result for arbitrary graded sequences of ideals.

Recall that if $\ba_{\bullet}$ and $\bb_{\bullet}$ are graded sequences
of ideals, then their intersection $\ba_{\bullet}\cap\bb_{\bullet}$
is defined by $\{\ba_m\cap\bb_m\}_m$. Similarly, their product
$\ba_{\bullet}\cdot\bb_{\bullet}$ is defined by
$\{\ba_m\cdot\bb_m\}_m$. It is clear that both
$\ba_{\bullet}\cap\bb_{\bullet}$
and $\ba_{\bullet}\cdot\bb_{\bullet}$ are graded sequences of ideals.
Moreover, if $\ba_{\bullet}$ and $\bb_{\bullet}$ are sequences of zero-dimensional
ideals, then so are $\ba_{\bullet}\cap\bb_{\bullet}$ and 
$\ba_{\bullet}\cdot\bb_{\bullet}$.

\begin{corollary}
If $\ba_{\bullet}$ and $\bb_{\bullet}$ are graded sequences
of zero-dimensional ideals, then
$$\vol(\ba_{\bullet}\cap\bb_{\bullet})^{1/n}\leq 
\vol(\ba_{\bullet}\cdot\bb_{\bullet})^{1/n}\leq \vol(\ba_{\bullet})^{1/n}
+\vol(\bb_{\bullet})^{1/n}.$$
\end{corollary}

\begin{proof}
Since $\ba_m\cdot\bb_m\subseteq\ba_m\cap\bb_m$ for all $m$,
we have a corresponding inequality between multiplicities. The
first inequality follows by dividing by $m^n$, taking the limit,
and applying Theorem~\ref{main}.

For the second one, by Theorem~\ref{main}, it is enough
to prove that $\ee(\ba_{\bullet}\cdot\bb_{\bullet})^{1/n}
\leq \ee(\ba_{\bullet})^{1/n}+ \ee(\bb_{\bullet})^{1/n}$.
For every $m$, Teissier's inequality (see \cite{teissier})
gives $\ee(\ba_m)^{1/n}\leq \ee(\ba_m)^{1/n}+\ee(\bb_m)^{1/n}$.
Dividing by $m$ and taking the limit, we get 
 our inequality.
\end{proof}

We show now the easy inequality $\vol(\ba_{\bullet})\leq
\ee(\ba_{\bullet})$ and use it to reduce the statement of Theorem~\ref{main}
to the case of a graded sequence of monomial ideals. The proof of that
case will be given in the next section.

\begin{lemma}\label{inequality}
If $\ba_{\bullet}$ is a graded sequence of zero-dimensional ideals,
then $\vol(\ba_{\bullet})\leq \ee(\ba_{\bullet})$.
\end{lemma}

\begin{proof}
It is enough to prove that for every $p$, we have
$$\limsup_{m\to\infty}\frac{n!\cdot \length(R/\ba_m)}{m^n}
\leq \frac{\ee(\ba_p)}{p^n}.$$
To see this, we show that for every integer
$k$, with $0\leq k<p$ we have
$$\limsup_{m\to\infty}\frac{n!\cdot
\length(R/\ba_{mp+k})}{(mp+k)^n}\leq\frac{\ee(\ba_p)}{p^n}.$$

Since $\ba_{\bullet}$ is a graded sequence, we have
$\ba_p^m\cdot\ba_k\subseteq\ba_{mp+k}$. 
If $\ba_p=R$, then 
this implies $\length(R/\ba_{mp+k})\leq \length(R/\ba_k)$ for every $m$,
hence
$$\lim_{m\to\infty}\frac{n!\cdot \length(R/\ba_{mp+k})}{(mp+k)^n}=0,$$
and we are done.

If $a_p\neq R$, then we can find $r\in\NN$ such that
$\ba_p^r\subseteq\ba_k$.
We deduce
$$\frac{n!\cdot \length(R/\ba_{mp+k})}{(mp+k)^n}
\leq\frac{n!\cdot \length(R/\ba_p^{m+r})}{(m+r)^n}\cdot\frac{(m+r)^n}{(mp+k)^n},$$
for every $m$. Since the right hand side has limit $\ee(\ba_p)/p^n$
when $m$ goes to infinity, we are done by taking $\limsup$ in the
above inequality.
\end{proof}

\begin{remark}
It follows from Theorem~\ref{main} that for every $p$,
$\vol(\ba_{\bullet})=\limsup_m\frac{n!\cdot \length(R/\ba_{mp})}{(mp)^n}$.
Indeed, this is  obvious since $\ee(\ba_{\bullet})$ is a limit.
If we assume that $\ba_p\subseteq \ba_q$ for $p>q$, then this
can be easily proved directly (see \cite{els}, Lemma~3.8).

 Once we know that $\vol(\ba_{\bullet})
=\limsup_m\frac{n!\cdot \length(R/\ba_{mp})}{(mp)^n}$,
the proof of the above lemma becomes even easier since 
$\ba_p^m\subseteq\ba_{pm}$ implies 
$$\frac{n!\cdot \length(R/\ba_{pm})}{(pm)^n}\leq \frac{n!\cdot \length(R/\ba_p^m)}
{p^nm^n}$$
for all $m$. Taking $\limsup$ with respect to $m$,
 we deduce $\vol(\ba_{\bullet})\leq
\ee(\ba_p)/p^n$.
\end{remark}

\begin{lemma}\label{reduction}
If Theorem~\ref{main} is known to be true for every graded sequence
 $\ba_{\bullet}$ of zero-dimensional monomial ideals
in a polynomial ring over a field, then the
theorem is true in general.
\end{lemma}

\begin{note*}
We have given all the definitions for a regular local ring $R$.
When we work with $R=K[X_1,\ldots,X_n]$, we refer to the corresponding
statements for the localization at $(X_1,\ldots,X_n)$.
However, since in this case all our ideals are supported at the origin,
this should cause no confusion, and we will simplify in this way the
notation. 
\end{note*}

\begin{proof}
Let $\hat{R}$ be the completion of $R$ at $\bm$.
If $\ba'_m=\ba_m\hat{R}$ for all $m$, it is clear that $\vol(\ba_{\bullet})
=\vol(\ba'_{\bullet})$ and $\ee(\ba_{\bullet})=\ee(\ba'_{\bullet})$.
Moreover, since $\ba_p\subseteq\bm^q$ if and only if $\ba'_p\subseteq
(\bm\hat{R})^q$, it is clear that Theorem~\ref{main} is true for
$\ba_{\bullet}$ if and only if it is true for $\ba'_{\bullet}$.

On the other hand, since $R$ is regular and contains a field, if
$K=R/\bm$, then $\hat{R}\simeq K[[X_1,\ldots,X_n]]$.
Reversing the previous argument, we see that it is enough to 
prove the theorem when $R=K[X_1,\ldots,X_n]$ and $\ba_p$
are ideals supported at the origin.

We consider now a deformation of $\ba_{\bullet}$ to a sequence of monomial
ideals. For example, pick a monomial order $>$ on $R$ and let
$\ba''_m={\rm in}_>(\ba_m)$ for all $m$
(see, for example, \cite{eisenbud}, Chapter 15, for initial monomial ideals).
 If $u={\rm in}_>(f)$
and $v={\rm in}_>(g)$ for $f\in\ba_p$ and $g\in\ba_q$, then 
$uv={\rm in}_>(fg)$ and $fg\in\ba_{p+q}$. Therefore $\ba''_{\bullet}$
is a graded sequence of monomial ideals, which are clearly supported
at the origin.

Moreover, we have $\length(R/\ba_p)=\length(R/\ba''_p)$ and $\ee(\ba''_p)\geq
\ee(\ba_p)$. The equality of lengths is well-known, 
while the inequality between
multiplicities can be seen as follows: since $(\ba_p'')^m\subseteq
 {\rm in}_>(\ba_p^m)$, 
we have $\length(R/(\ba_p'')^m)\geq \length(R/\ba_p^m)$
for all $m$. Dividing by $m^n$ and taking the limit with respect to $m$
gives the inequality.

We deduce $\vol(\ba_{\bullet})=
\vol(\ba''_{\bullet})$ and $\ee(\ba''_{\bullet})
\geq \ee(\ba_{\bullet})$. Using Lemma~\ref{inequality}, we have
$$\vol(\ba''_{\bullet})=\vol(\ba_{\bullet})\leq \ee(\ba_{\bullet})
\leq \ee(\ba''_{\bullet}).$$ 
Since assertion $(1)$ in the theorem is true for 
$\ba''_{\bullet}$, we deduce this assertion  
for $\ba_{\bullet}$.

For the proof of $(2)$, note that one implication is trivial.
Namely, if $\ba_p\subseteq\bm^{[p/q]}$ for all $p$, then
$$\length(R/\ba_p)\geq \length(R/\bm^{[p/q]})={{[p/q]+n-1}\choose{n}}.$$
Dividing by $p^n$ and taking $\limsup$ gives
$\vol(\ba_{\bullet})\geq (1/q)^n>0$.

For the converse, once we know the theorem for $\ba_{\bullet}''$,
it is enough to show that we can make the deformation from $\ba_{\bullet}$
to $\ba_{\bullet}''$ such that for every $p$
 and $r$, $\ba_p\subseteq\bm^r$ if
$\ba''_p\subseteq\bm^r$. This is clear if $\ba_m$ is homogeneous for every
$m$.
In the general case, consider the graded sequence of ideals
$\widetilde{\ba}_m=(l(f)\mid f\in\ba_m)$ where 
$l(f)$ is the sum of the terms in $f$ of smallest degree. It is clear that
$\widetilde{\ba}_{\bullet}$ is a graded system of homogeneous ideals such
that
$\vol(\ba_{\bullet})=\vol(\widetilde{\ba}_{\bullet})$ and 
$\ba^p\subseteq\bm^r$
if and only if $\widetilde{\ba}_p\subseteq\bm^r$. Since we know $(2)$
for $\widetilde{\ba}_{\bullet}$, we deduce it for $\ba_{\bullet}$, and
this completes the proof of the lemma. 
\end{proof}

It follows from the above proof that the computation of
$\ee(\ba_{\bullet})$ can be reduced to the case of a graded sequence of
monomial ideals. We state this as a separate corollary.

\begin{corollary}\label{reduction_of_computation}
Let $\ba_{\bullet}$ be a graded sequence of monomial ideals 
in $R=K[X_1,\ldots,X_n]$, supported at the origin. If $>$ is a monomial
order on the monomials in $R$, and if $\ba'_m=\lin_{>}(\ba_m)$ for all $m$,
then $\ee(\ba_{\bullet})=\ee(\ba'_{\bullet})$. In particular, for every
ideal $I\subset R$ supported at the origin, we have
$$\ee(I)=\lim_{m\to\infty}\frac{\ee(\lin_{>}(I^m))}{m^n}.$$
\end{corollary}

\section{Graded sequences of monomial ideals}

In this section we finish the proof of Theorem~\ref{main},
by proving it for graded sequences of monomial ideals.
In this case, we prove a stronger statement involving the asymptotic
multiplier ideals of $\ba_{\bullet}$. More precisely, we prove that
the full conclusion of Proposition~3.11 in \cite{els} remains
true for an arbitrary graded sequence of zero-dimensional monomial ideals.

Note that since we work over a field of arbitrary characteristic,
the usual results concerning multiplier ideals do not apply.
 On the other hand, since in this
section we are concerned only with monomial ideals, the characteristic
of the field does not play any role and we could always reduce ourselves
to a field of characteristic zero. However, in order to underline
the elementary nature of the arguments, we will define directly in this
case multiplier ideals and deduce the basic property that we need,
 the subadditivity, directly from definition.
We start with some general considerations. Recall that we work in a ring 
$R$ which is either a regular local ring or a polynomial ring over a field.

\begin{definition}
A reverse-graded sequence
of ideals is a family of
ideals $\bb_{\bullet}=\{\bb_m\}_{m\in\NN^*}$ such that
\begin{enumerate}
\item If $p>q$, then $\bb_p\subseteq\bb_q$.
\item $\bb_{p+q}\subseteq\bb_p\cdot\bb_q$, for every $p$, $q\in\NN^*$.
\end{enumerate}
If $\ba_{\bullet}$ is a graded sequence of ideals, then we say that
$\bb_{\bullet}$ dominates $\ba_{\bullet}$ if $\ba_m\subseteq\bb_m$ for
every $m\in\NN^*$. 
\end{definition}

We have the following lemma, which plays an analogous role with
Lemma~\ref{limit}.

\begin{lemma}\label{limit1}
If $\{\beta_m\}_{m\in\NN^*}$ is a non-decreasing sequence of non-negative
real numbers, such that $\beta_{mp}\geq m\beta_p$ for all $m$ and $p$,
then $\lim_{m\to\infty}\beta_m/m=\sup_m\beta_m/m$.
\end{lemma}

\begin{proof}
Let $M=\sup_m\beta_m/m$. Suppose that $M<\infty$, the case
$M=\infty$ being analogous.
For every $\epsilon>0$, pick $p$ such that $\beta_p/p\geq M-\epsilon/2$.
It is enough to show that for every integer $q$, with $0\leq q<p$, we have
$\beta_{mp+q}/(mp+q)\geq M-\epsilon$ for $m\gg 0$.
Since we have
$$\frac{\beta_{mp+q}}{mp+q}\geq\frac{\beta_{mp}}{mp+q}
\geq\frac{m\beta_p}{mp+q}\geq
(M-\epsilon/2)\cdot\frac{mp}{mp+q},$$
and since the right hand side goes to $M-\epsilon/2$ when $m$ goes to
infinity,
it follows that for $m\gg 0$ we have $\beta_{mp+q}/(mp+q)^n\geq
M-\epsilon$.
\end{proof}

\begin{corollary}
If $\bb_{\bullet}$ is a reverse-graded sequence of zero-dimensional
ideals, then $\lim_{m\to\infty}\ee(\bb_m)/m^n=\sup_m\ee(\bb_m)/m^n$.
\end{corollary}

\begin{proof}
Apply the above lemma to the sequence $\beta_m=\ee(\bb_m)^{1/n}$, noting
that 
$\bb_{mp}\subseteq\bb_p^m$ implies $\ee(\bb_{mp})\geq\ee(\bb_p)\cdot m^n$.
\end{proof}

\begin{definition}
If $\bb_{\bullet}$ is a reverse-graded sequence of zero-dimensional ideals,
then the volume of $\bb_{\bullet}$ is defined by
the same formula as $\vol(\ba_{\bullet})$:
$$\vol(\bb_{\bullet}):=\limsup_m\frac{n!\cdot \length(R/\bb_m)}{m^n}.$$
The multiplicity of $\bb_{\bullet}$ is defined by
$$\ee(\bb_{\bullet}):=\lim_{m\to\infty}\frac{\ee(b_m)}{m^n}.$$
\end{definition}

\begin{lemma}\label{one_equality}
If $\bb_{\bullet}$ is a reverse-graded sequence of zero-dimensional ideals,
then we have $\vol(\bb_{\bullet})=\limsup_{m\to\infty}
\frac{n!\cdot \length(R/\bb_{mp})}{(mp)^n}$, for every $p\in\NN^*$.
\end{lemma}

\begin{proof}
It is clear that $\vol(\bb_{\bullet})\geq L:=\limsup_m\frac{n!\cdot
\length(R/\bb_{mp})}
{(mp)^n}$. For the reverse inequality we use again the standard argument:
for $\epsilon>0$, let $m_0$ be such that $\frac{n!\cdot \length(R/\bb_{mp})}
{(mp)^n}\leq L+\epsilon/2$ for all $m\geq m_0$. It is enough to prove that
for every integer $q$, with $0\leq q< p$, we have
$\frac{n!\cdot \length(R/\bb_{mp+q})}{(mp+q)^n}\leq L+\epsilon$ for all $m\gg
0$.
Since $\bb_{(m+1)p}\subseteq\bb_{mp+q}$, we have
$$\frac{n!\cdot \length(R/\bb_{mp+q})}{(mp+q)^n}\leq \frac{n!\cdot
\length(R/\bb_{(m+1)p})}
{(m+1)^np^n}\cdot\frac{(mp+p)^n}{(mp+q)^n}.$$
We are done if $m\geq\max\{m_0,m_1\}$, where $m_1$ is such that
$(mp+p)^n/(mp+q)^n\leq (L+\epsilon)/(L+\epsilon/2)$ for $m\geq m_1$.
\end{proof}

\begin{lemma}\label{three_ineq}
Let $\ba_{\bullet}$ be a graded sequence
of zero-dimensional ideals and let $\bb_{\bullet}$ be a reverse-graded
sequence dominating $\ba_{\bullet}$.
We have the following inequalities: 
$$\ee(\bb_{\bullet})\leq \vol(\bb_{\bullet})\leq \vol(\ba_{\bullet})\leq 
\ee(\ba_{\bullet}).$$
\end{lemma}

\begin{proof}
Since we have proved in Lemma~\ref{inequality}
that $\vol(\ba_{\bullet})\leq \ee(\ba_{\bullet})$ and since $\vol(\bb_{\bullet})
\leq \vol(\ba_{\bullet})$ follows trivially from $\ba_m\subseteq\bb_m$ for
all
$m$, it is enough to show that $\ee(\bb_{\bullet})\leq \vol(\bb_{\bullet})$.

Fix $p$. Since $\bb_{mp}\subseteq \bb_p^m$ for all $m$,
we deduce $\length(R/\bb_{mp})\geq \length(R/\bb_p^m)$. 
Multiplying by $n!/(mp)^n$ and taking $\limsup$ when
$m$ goes to infinity, we deduce by Lemma~\ref{one_equality}
 that $\vol(\bb_{\bullet})
\geq \ee(\bb_p)/p^n$. Since this holds for every $p$, we get
$\vol(\bb_{\bullet})\geq\sup_p\ee(\bb_p)/p^n=\ee(\bb_{\bullet})$.
\end{proof}

\begin{remark}\label{zero}
Since $\ee(\bb_{\bullet})=\sup_m\ee(\bb_m)/b^m$, it follows
that $\ee(\bb_{\bullet})=0$ if and only if $\bb_m=R$ for every $m$.
\end{remark}

\smallskip

To check equality in Lemma~\ref{three_ineq} we will use the following
criterion from \cite{els}. Suppose that $\ba_{\bullet}$
is a graded sequence of zero-dimensional ideals and $\bb_{\bullet}$
is a reverse-graded sequence which dominates $\ba_{\bullet}$.
It is easy to check that the set of colon ideals $\{\ba_m:\bb_m\}_m$
forms a graded sequence of ideals, which we denote by 
$\ba_{\bullet}:\bb_{\bullet}$. Note that since $\ba_m$ is zero-dimensional,
so is $\ba_m:\bb_m$.

\begin{lemma}\label{key}{\rm (\cite{els}, 3.11)}
With the above notation, if $\ee(\ba_{\bullet}:\bb_{\bullet})=0$,
then 
$$\ee(\bb_{\bullet})=\vol(\bb_{\bullet})=
\vol(\ba_{\bullet})=\ee(\ba_{\bullet}).$$
\end{lemma}

\begin{proof}
We recall the proof for completeness.
By Lemma~\ref{three_ineq}, it is enough to prove that $\ee(\ba_{\bullet})
\leq \ee(\bb_{\bullet})$. Let $\bc_m=\ba_m:\bb_m$. Since we have
$\bb_m\cdot\bc_m\subseteq\ba_m$ for all $m$, using Teissier's
inequality \cite{teissier}, we deduce
$$\ee(\ba_m)^{1/n}\leq \ee(\bb_m\cdot\bc_m)^{1/n}\leq \ee(\bb_m)^{1/n}
+\ee(\bc_m)^{1/n}.$$
If we divide by $m$ and take the limit when $m$ goes to infinity,
the hypothesis implies $\ee(\ba_{\bullet})\leq \ee(\bb_{\bullet})$.
\end{proof}

\bigskip

As in \cite{els},
the sequence $\bb_{\bullet}$ we use is given by the
asymptotic multiplier ideals of $\ba_{\bullet}$. From now on we
fix a graded sequence $\ba_{\bullet}$
consisting of monomial ideals in $R=K[X_1,\ldots,X_n]$,
which are supported at the origin. If $u=(u_i)_i\in\NN^n$, we use the
notation $X^u=\prod_iX_i^{u_i}$.

\begin{definition}
Let $\ba\subseteq R=K[X_1,\ldots,X_n]$ be a monomial ideal and
$P_{\ba}$ its Newton polyhedron, i.e., $P_{\ba}$ is the convex hull
of $\{u\in\NN^n\vert X^u\in\ba\}$. If $\lambda\in\QQ_+^*$, then
the multiplier ideal of $\ba$ with coefficient $\lambda$ is the monomial
ideal
$$\cI(\lambda\cdot\ba):=(X^u\vert u\in\NN^n, u+e\in{\rm Int}(\lambda \cdot
P_{\ba})),$$
where $e=(1,\ldots,1)\in\NN^n$.
\end{definition}

\begin{remark}
The usual definition of multiplier ideals is different (see \S 3),
 and it is a
theorem
of Howald from \cite{howald} that for a monomial ideal we have this
expression.
\end{remark}

Suppose now that $\ba_{\bullet}$ is a graded sequence of monomial ideals 
in $R$. It is clear that if $\lambda\in\QQ_+^*$, and if $p$, $q\in\NN^*$,
then 
$$\cI(\lambda/p\cdot\ba_p)\subseteq\cI(\lambda/pq\cdot\ba_{pq}).$$
Indeed, this follows since $P_{\ba_p}\subseteq (1/q)P_{\ba_q}$,
as $qP_{\ba_p}\subseteq P_{\ba_p^q}\subseteq P_{\ba_{pq}}$.
It is clear from this that the set $\{\cI(\lambda/p\cdot\ba_p)\}_p$
has a unique maximal element, called the asymptotic multiplier ideal of
$\ba_{\bullet}$ with coefficient $\lambda$, and
 denoted by $\cI(\lambda\cdot\parallel\ba_{\bullet}\parallel)$.

Given the graded sequence $\ba_{\bullet}$, we take $\bb_m=\cI(m\cdot
\parallel \ba_{\bullet}\parallel)$. We then have the following

\begin{lemma}
Let $\ba_{\bullet}$ be a graded sequence of monomial ideals. 
With the above definition, $\bb_{\bullet}$ is a reverse-graded sequence
of ideals dominating $\ba_{\bullet}$.
\end{lemma}

\begin{proof}
It follows from definition that $\ba_m\subseteq\cI(\ba_m)\subseteq\bb_m$.
Moreover, 
it is clear that if $\lambda<\mu$, then $\cI(\mu\cdot\ba)
\subseteq\cI(\lambda\cdot\ba)$ for every $\ba$. This immediately
implies $\bb_q\subseteq\bb_p$ for $p<q$.

The last property we need for $\bb_{\bullet}$ follows from the general
subadditivity theorem (see \cite{del}). In the case of monomial
ideals it is very easy to give a direct proof. Note that it is a formal
consequence of the following assertion: 
if $\ba$ and $\ba'$ are monomial ideals, and if
$\lambda\in\QQ_+^*$, we have 
$$\cI(\lambda\cdot(\ba\cdot\ba'))
\subseteq\cI(\lambda\cdot\ba)\cdot\cI(\lambda\cdot\ba').$$

In order to prove this, suppose that
$X^u\in\cI(\lambda\cdot(\ba\cdot\ba'))$, i.e.,
$u+e\in{\rm Int}(\lambda\cdot(P_{\ba}+P_{\ba'}))$.
This means that we can write $u+e=\lambda(v+w)$, where we may assume,
for example, that $v\in{\rm Int}(P_{\ba})$ and $w\in P_{\ba'}$.

For $x\in\RR$, denote by $\{x\}$ the smallest integer $m$ such that $m>x$.
Note that $\{x\}\leq x+1$. If $\alpha=(\alpha_i)_i\in\RR^n$, we put 
$\{\alpha\}=(\{\alpha_i\})_i$. 

Take $w'=\{\lambda w\}-e$ and $v'=u-w'$. Therefore we have $u=v'+w'$.
By definition, we have $X^{w'}\in\cI(\lambda\cdot P_{\ba'})$.
Moreover, we have $X^{v'}\in\cI(\lambda\cdot P_{\ba})$. Indeed,
$e+v'=(\lambda w+e-\{\lambda w\})+\lambda v$, and the first term is in
$\RR_+^n$, while $\lambda v\in{\rm Int}(\lambda P_{\ba})$. This completes the
proof
of the lemma.
\end{proof}

\smallskip

The following is the main result of this section.

\begin{theorem}\label{all_equal}
Let $\ba_{\bullet}$ be a graded sequence of 
zero-dimensional monomial ideals.
 If $\bb_{\bullet}$ is the corresponding
sequence of asymptotic multiplier ideals, then
$$\ee(\bb_{\bullet})= \vol(\bb_{\bullet})=\vol(\ba_{\bullet})
=\ee(\ba_{\bullet}).$$
\end{theorem}

Granted this, we can finish the proof of the result we have stated
in the previous section.

\begin{proof}[Proof of Theorem~\ref{main}]
By Lemma~\ref{reduction}, we may assume that
all $\ba_m$ are monomial ideals in $R=K[X_1,\ldots,X_n]$.
The assertion in $(1)$ follows from the more precise statement
in Theorem~\ref{all_equal} above. Moreover, we have seen
in the proof of Lemma~\ref{reduction} that the only nontrivial
implication in $(2)$ is that if $\ee(\ba_{\bullet})>0$, then there is $q\in\NN^*$
such that $\ba_p\subseteq\bm^{[p/q]}$ for all $p$.

By Theorem~\ref{all_equal} above, $\ee(\ba_{\bullet})>0$
implies $\ee(\bb_{\bullet})>0$, i.e., there is $q$ such that
$\bb_q\subseteq\bm=(X_1,\ldots,X_n)$. Since this implies
$$\ba_p\subseteq\bb_p\subseteq\bb_{q[p/q]}\subseteq\bb_q^{[p/q]}
\subseteq\bm^{[p/q]},$$
we are done.    
\end{proof}

\bigskip

Before giving the proof of Theorem~\ref{all_equal}, we need some preparation.
We start by interpreting $\ee(\ba_{\bullet})$ in terms
of the polyhedra involved. If $\ba_{\bullet}$ is a graded sequence
of zero-dimensional monomial ideals,
let $Q_m$ be the closure of
$\RR_+^n\setminus P_{\ba_m}$
(if $\ba_m=R$, then we take $Q_m=\{0\}$).
 It is clear that $Q_m$ is compact for every
$m$. Moreover, the condition that $\ba_{\bullet}$ is a graded sequence
of ideals implies 
\begin{equation}\label{eq1}
Q_{p+q}\subseteq Q_p+ Q_q,
\end{equation}
for every $p$ and $q$. In particular, $(1/p)Q_p\subseteq (1/q)Q_q$
if $q$ divides $p$.

Indeed, if 
$u\in \RR_+^n\setminus P_{\ba_{p+q}}$, and if 
$v\in P_{\ba_p}\cap Q_p$
is such that $u-v\in\RR_+^n$, then
$u-v\in Q_q$. Note that we can choose such $v$,
unless $u\in Q_p$, in which case we have $u\in Q_p+Q_q$ trivially.
We deduce now equation~(\ref{eq1}), since the right hand side is closed. 

Let $Q:=\bigcap_{m\in\NN^*}(1/m)Q_m$. It is clear that $Q$ is compact.
Recall the well-known fact that $\ee(\ba_m)=n!\vol(Q_m)$. The following
lemma implies the analogous equality for a graded sequence:
$\ee(\ba_{\bullet})=n!\vol(Q)$.

\begin{lemma}\label{interpretation}
For every neighbourhood $U$ of $Q$, we have $(1/m)Q_m\subseteq U$
for $m\gg 0$. In particular,
$\vol(Q)=\lim_{m\to\infty}\frac{\vol(Q_m)}{m^n}$,
hence $\ee(\ba)=n!\vol(Q)$.
\end{lemma}

\begin{proof}
Fix an open neighbourhood $W$ of $Q$ such that $\overline{W}$
is compact and contained in $U$. Moreover, since $\lambda Q\subseteq Q$
for every $\lambda$ with $0\leq\lambda\leq 1$, we may assume
that $W$ also has this property.

Since all $Q_m$ are closed and lie in a bounded domain,
we can find $m_1,\ldots,m_k\in\NN$ such that $\bigcap_{1\leq i\leq k}
(1/m_i)Q_{m_i}\subseteq W$. If we pick $m_0$ such that $m_0$
is divisible by $m_i$
for $1\leq i\leq k$, it follows that $(1/m)Q_m\subseteq W$ if
$m_0$ divides $m$. In order to finish, it is enough to show that
for every integer 
$q$, with $0<q<m_0$, we have $(1/lm_0+q)Q_{lm_0+q}\subseteq U$
for $l\gg 0$.

Let $U_0$ be an open neighbourhood of $0$ such that $\overline{W}+U_0
\subseteq U$. If we choose $\mu>0$ such that $\mu\cdot(1/q)Q_q\subseteq
U_0$
and if $l_0$ is such that $q/(l_0m_0+q)<\mu$, then it follows from
the inclusion~$(\ref{eq1})$
 and our conditions on $W$, $U_0$, $l_0$ and $m_0$ that
$(1/lm_0+q)Q_{lm_0+q}\subseteq U$ for all $l\geq l_0$.
\end{proof}

For the proof of Theorem~\ref{all_equal} we will use 
Lemma~\ref{key}. Note however that, as the following example
shows, an arbitrary graded sequence of monomial ideals does
not satisfy the hypothesis of that lemma.

\begin{example}
Let $\ba_{\bullet}$ be the graded sequence 
in $R=K[x,y]$ defined by
$\ba_m=(x^m,y^m)$. It is easy to see that $\bb_m=(x,y)^{m-1}$.
Since $(xy)^p\in\bb_{2p+1}$, 
it follows that $(\ba_{2p+1}:\bb_{2p+1})\subseteq
 (x^{p+1},y^{p+1})$. Therefore $\ee(\ba_{2p+1}:\bb_{2p+1})\geq (p+1)^2$ 
for all
$p$, hence $\ee(\ba_{\bullet}:\bb_{\bullet})\geq 1/4$.
\end{example}

\smallskip

We show now how to associate to a graded sequence of 
monomial ideals $\ba_{\bullet}$
another closely related such sequence $\ba'_{\bullet}$,
which satisfies the hypothesis of Lemma~\ref{key}.

If $\ba_{\bullet}$ is a graded sequence of monomial ideals and if
$m$ is fixed, then consider the family of ideals 
$\{\ba'_{m,r}\}_r$, where 
$\ba'_{m,r}=(X^u\vert u\in \NN^n\cap (1/r)P_{\ba_{mr}})$. It is clear that
if $r$ divides $p$, then $\ba'_{m,r}\subseteq\ba'_{m,p}$. Since $R$
is Noetherian, it follows that there is a unique maximal element among
$\{\ba'_{m,r}\}_r$, which we denote by $\ba'_m$. It is clear that
$\ba'_{\bullet}$ is a graded sequence of monomial ideals 
such that $\ba_m\subseteq\ba'_m$ for all $m$. 

\begin{lemma}\label{same_multiplier}
If $\ba_{\bullet}$ and $\ba'_{\bullet}$ are as above, and if
$\bb_{\bullet}$ and $\bb'_{\bullet}$ are the corresponding sequences
of asymptotic multiplier ideals, then $\bb_{\bullet}=\bb'_{\bullet}$.
\end{lemma}

\begin{proof}
 Since $\ba_m\subseteq\ba'_m$ for every $m$, it is clear that
$\bb_m\subseteq\bb'_m$ for every $m$. For the reverse inclusion,
we have to prove that $\cI(1/p\cdot\ba'_{pm})\subseteq\,\,\,\,\,\,\,\,
\cI(m\cdot\parallel
\ba_{\bullet}\parallel)$ for every $p$ and $m$. Suppose that $q$ is such that 
$\ba'_{pm}=\ba'_{pm,q}$. If $X^u\in\cI(1/p\cdot\ba'_{pm})$, then
$u+e\in(1/p){\rm Int}(P_{\ba'_{pm}})$, hence $pq(u+e)\in{\rm
Int}(P_{\ba_{pmq}})$.
It follows from definition that $X^u\in\cI(1/pq\cdot\ba_{pmq})\subseteq
\cI(m\cdot\parallel\ba_{\bullet}\parallel)$.
\end{proof}

\begin{lemma}\label{same_multiplicity}
If $\ba_{\bullet}$ is a graded sequence of zero-dimensional
monomial ideals, then
we have, with the above notation, $\ee(\ba_{\bullet})=\ee(\ba'_{\bullet})$.
\end{lemma}

\begin{proof}
We use Lemma~\ref{interpretation}. With the notation in that lemma,
it is enough to prove that if $Q$ and $Q'$ are the compact sets
corresponding
to these two graded sequences of ideals, we have $Q\cap(\RR^*_+)^n
=Q'\cap(\RR^*_+)^n$. Since
$\ba_m\subseteq\ba'_m$ for every $m$, we clearly have $Q'\subseteq Q$. We
show now that $Q\cap(\RR_+^*)^n\subseteq Q'\cap(\RR^*_+)^n$.
 Suppose that $u\in Q\cap(\RR^*_+)^n$, but $u\not\in Q'$.
Then there is $m$ such that $u\in(1/m){\rm Int}(P_{\ba'_m})$. Let $p$
be such that $\ba'_m=\ba'_{m,p}$. Since $p\cdot P_{\ba'_{m,p}}\subseteq
P_{\ba_{pm}}$,
we deduce 
$u\in(1/mp){\rm Int}(P_{\ba_{mp}})\subseteq\RR^n\setminus Q$, a
contradiction.
\end{proof}

\begin{definition}
We say that a graded sequence of monomial ideals $\ba_{\bullet}$
is saturated if $\ba_{\bullet}=\ba_{\bullet}'$.
\end{definition}

\begin{lemma}\label{sat}
If $\ba_{\bullet}$ is a graded sequence of monomial ideals, then 
$\ba_{\bullet}'$ is saturated.
\end{lemma}

\begin{proof}
We have to prove that for every $p$ and $m$, if $u\in\NN^n$
is such that $pu\in P_{\ba'_{mp}}$, then $X^u\in\ba'_m$.
Let $r$ be such that $\ba'_{mp}=\ba'_{mp,r}$. By definition,
$pu\in P_{\ba'_{mp,r}}$ implies $rpu\in P_{\ba_{mpr}}$, hence
$X^u\in\ba'_{m,pr}\subseteq\ba'_m$.
\end{proof}

\begin{lemma}\label{sat_is_tight}
If $\ba_{\bullet}$ is a saturated graded sequence of
zero-dimensional monomial ideals, and if $\bb_{\bullet}$
is the corresponding sequence of asymptotic multiplier ideals, then
we have $\ee(\ba_{\bullet}:\bb_{\bullet})=0$.
\end{lemma}

\begin{proof}
Let $\bc_m=(\ba_m:\bb_m)$, and we first show that $X^e\in\bigcap_m
\bc_m$. If $X^u\in\bb_m$, and if $p$ is such that 
$\bb_m=\cI((1/p)\cdot\ba_{pm})$, we have in particular $u+e\in
(1/p)P_{\ba_{pm}}$.
Therefore $X^e\cdot X^u\in\ba'_{m,p}\subseteq\ba'_m=\ba_m$, since
$\ba_{\bullet}$ is saturated. Hence $X^e\in\bc_m$.

It is now easy to see that $\ee(\bc_{\bullet})=0$. Indeed, 
let us consider for every $i$, the polynomial ring
$R_i=K[X_1,\ldots,\hat{X_i},\ldots,X_n]$ for every $i$,
and let ${\bc}_{m,i}=\bc_m\cap R_i$. It is clear that
$\bc_{{\bullet},i}$ is a graded sequence of monomial ideals in $R_i$.
Moreover, there is a constant $C$ depending only on $n$ such that
$$\ee(R/\bc_m)\leq C\left(\sum_{i=1}^n\ee(R_i/\bc_{m,i})+1\right),$$ 
for every $m$.
Dividing by $m^n$ and taking the limit when $m$ goes to infinity gives
$\ee(\bc_{\bullet})=0$, since $\dim(R_i)=n-1$ for every $i$.
\end{proof}

We can give now the proof of
Theorem~\ref{all_equal}

\begin{proof}[Proof of Theorem~\ref{all_equal}]
By Lemma~\ref{three_ineq}, it is enough to prove that $\ee(\ba_{\bullet})
=\ee(\bb_{\bullet})$. Using Lemmas~\ref{same_multiplier} and
\ref{same_multiplicity},
we see that it is enough to prove that
$\ee(\ba'_{\bullet})=\ee(\bb_{\bullet}')$, where $\bb_{\bullet}'$
is the sequence of asymptotic multiplier ideals corresponding to 
$\ba'_{\bullet}$.
Therefore, by Lemma~\ref{sat}, we may assume that $\ba$ is saturated. 
Lemma~\ref{sat_is_tight} shows that $\ba_{\bullet}$ satisfies
the hypothesis of
Lemma~\ref{key}, so we are done.
\end{proof}

\begin{question}
A basic question is whether the assertion in Theorem~\ref{all_equal}
remains true for arbitrary graded sequences of zero-dimensional ideals 
(assuming that we are in a setting where
we have available the theory of multiplier ideals). 
We will see in Theorem~\ref{equality_for_lc}
that the analogous assertion is true if we replace the multiplicity by
the log canonical threshold: we have $\lc(\ba_{\bullet})=\lc(\bb_{\bullet})$.
\end{question}

\section{The log canonical threshold of a graded system of ideals}

We apply now the ideas used in the previous sections to the
study of log canonical thresholds. We suppose that we are in a geometric
situation: let $X$ be a smooth variety over an algebraically closed
field $k$ of characteristic zero, and let $R$ be the local ring of $X$
at a (not necessarily closed) point.

We recall briefly the definition of multiplier ideals, and refer 
for details and basic properties to \cite{lazarsfeld}.
Let $\ba\subseteq R$ be a non-zero ideal and $V(\ba)\subseteq
X=\Spec\,R$, the subscheme defined by $\ba$.
Let $f:X'\longrightarrow
X$ be a log resolution for $(X, V(\ba))$, i.e., 
a proper, birational morphism, with $X'$ smooth, and
such that $f^{-1}(V(\ba))\cup {\rm Ex}(f)$ is a divisor with simple normal
crossings (${\rm Ex}(f)$ denotes the exceptional locus of $f$). 
Let $K_{X'/X}$ be the relative canonical divisor of $f$.

If $\lambda\in\QQ_+^*$, and if $D=[\lambda\cdot f^{-1}(V(\ba))]$,
then the multiplier ideal of $\ba$ with coefficient
$\lambda$ is 
$$\cI(\lambda\cdot\ba):=f_*(\cO_{X'}(K_{X'/X}
-D)).$$
One shows that the definition does not depend on the particular resolution,
and this fact can be conveniently expressed as follows. Suppose that
$E$ is a divisor with center on $X$, i.e., it is a divisor on some smooth
model $\widetilde{X}$ over $X$. We identify $E$ with the corresponding
discrete valuation ring $\cO_{X',E}$ and ${\rm ord}_E$ will
denote the induced valuation. If $\ba'$ is an ideal in $R$, then
we put ${\rm ord}_E(\ba'):=\inf\{{\rm ord}_E(u)\mid u\in\ba'\}$. 
With this notation, if $u\in R$, then
$u\in\cI(\lambda\cdot\ba)$ if and only if for every $E$ as above, we have
$${\rm ord}_E(u)>{\rm ord}_E(\ba)-{\rm ord}_E(K_{X'/X})-1.$$

Going from multiplier ideals to asymptotic multiplier ideals involves
the same process as the one we sketched in the previous section
(see \cite{lazarsfeld} for details).
As before, we put $\bb_m=\cI(m\cdot\parallel\ba_{\bullet}\parallel)$.
It follows from the Subadditivity Theorem (see \cite{del}) that
$\bb_{\bullet}$ is a reverse-graded sequence of ideals dominating 
$\ba_{\bullet}$.

For a non-zero ideal $\ba\subseteq R$, 
we denote by $\lc(\ba)$ the log canonical
threshold of the subsubscheme $V(\ba)$ 
(see \cite{kollar} for basic facts
about log canonical thresholds). It is defined as follows. 
If $f$ is a log resolution for $(X, V(\ba))$, as above, 
we write
$f^{-1}(V(\ba))=\sum_i\alpha_iE_i$ and 
 $K_{X'/X}=\sum_i\gamma_iE_i$, and then 
$$\lc(\ba):=\inf_i\frac{\gamma_i+1}{\alpha_i}.$$
In terms of multiplier ideals, we have 
$$\lc(\ba)=\sup\{\lambda>0\mid\cI(\lambda\cdot\ba)=R\}.$$
Note that $\lc(\ba)\in\QQ_+^*$, for every non-zero ideal $\ba$.

Recall the characterization of multiplier for monomial ideals, due
to Howald, which we have used in the previous section. It follows
from that description that if $\ba$ is a monomial ideal with
Newton polyhedron $P_{\ba}$, and if $e=(1,\ldots,1)$, then
$$1/\lc(\ba)=\inf\{\mu>0\mid\mu\cdot e\in P_{\ba}\}.$$

\begin{lemma}\label{lc_ineq}
If $\ba$ is a graded sequence of ideals, then
for every $p$ and $q$, we have $\frac{1}{\lc(\ba_{p+q})}\leq
\frac{1}{\lc(\ba_p)}+\frac{1}{\lc(\ba_q)}$.
\end{lemma}

\begin{proof}
Since $\ba_p\cdot\ba_q\subseteq\ba_{p+q}$, we deduce
$1/\lc(\ba_{p+q})\leq 1/\lc(\ba_p\cdot\ba_q)$. The statement of the lemma
follows once we show that for arbitrary ideals $\ba$ and $\bb$,
we have the following analogue of Teissier's inequality:
$$\frac{1}{\lc(\ba\cdot\bb)}\leq\frac{1}{\lc(\ba)}+\frac{1}{\lc(\bb)}.$$
Indeed, suppose that $f:X'\longrightarrow X=\Spec(R)$ is a log resolution
for
$(X,V(\ba)\cup V(\bb))$. If we write $f^{-1}(V(\ba))=\sum_i\alpha_iE_i$,
$f^{-1}(V(\bb))=\sum_i\beta_iE_i$ and $K_{X'/X}=\sum_i\gamma_iE_i$, then
$f^{-1}(\ba\bb)=\sum_i(\alpha_i+\beta_i)E_i$, and
$$\sup_i\frac{\alpha_i+\beta_i}{\gamma_i+1}\leq
\sup_i\frac{\alpha_i}{\gamma_i+1}+\sup_i\frac{\beta_i}{\gamma_i+1},$$
which is precisely our assertion.
\end{proof}

\begin{definition}
If $\ba_{\bullet}$ is a graded sequence of ideals in $R$, we define
the log canonical threshold of $\ba_{\bullet}$ by
$\lc(\ba_{\bullet}):=\lim_{m\to\infty}m\cdot \lc(\ba_m)$. 
By Lemma~\ref{lc_ineq}, we may apply Lemma~\ref{limit}
to the sequence $\{1/\lc(\ba_m)\}_m$ to see that $\lc(\ba_{\bullet})$
exists in $\RR_+^*\cup\{\infty\}$, and it is equal to
$\sup\{m\cdot \lc(\ba_m)\mid m\in\NN^*\}$.
\end{definition}

\begin{remark}
If $\ba_{\bullet}$ is a graded sequence of ideals as above, 
then
$$\lc(\ba_{\bullet})=\sup\{\lambda\in\QQ_+^*
\vert\cI(\mu\cdot\parallel\ba_{\bullet}
\parallel)=R\,\text{for all}\,\mu<\lambda\}.$$ 
Indeed, we have $\cI(\mu\cdot\parallel \ba_{\bullet}\parallel)=R$
if and only if there is some $p$ such that $\cI(\mu/p\cdot\ba_p)=R$,
which means that $\mu<p\cdot \lc(\ba_p)$.
\end{remark}

\begin{remark}
It follows from the above remark that
if $\ba_{\bullet}$ is a graded sequence of ideals in $R$,
and if $\bb_{\bullet}$ is the corresponding sequence of asymptotic
multiplier ideals,
then $\lc(\ba_{\bullet})=\infty$ if and only if $\bb_m=R$ for all $m$.
Recall that if $\ba_m$ 
is zero-dimensional for every $m$, then Theorem~\ref{main} 
shows that this is the case if and only if $\ee(\ba_{\bullet})=0$.
\end{remark}

\begin{definition}
If $\ba_{\bullet}$ is a graded sequence of ideals in $R$, and if
$\bb_{\bullet}$ is the corresponding sequence of asymptotic multiplier
ideals, then we define the log canonical threshold of $\bb_{\bullet}$
by $\lc(\bb_{\bullet}):=\lim_{m\to\infty}m\cdot \lc(\bb_m)$. 
It follows from Lemma~\ref{limit1} applied for $\beta_m=1/\lc(\bb_m)$
that $\lc(\bb_{\bullet})$
exists in $\RR_+\cup\{\infty\}$ and it is equal to
$\inf\{m\cdot \lc(\bb_m)\mid m\in\NN^*\}$.
\end{definition}

The following result shows that with respect to the log canonical threshold,
the sequences $\ba_{\bullet}$ and $\bb_{\bullet}$ grow in the same way.

\begin{theorem}\label{equality_for_lc}
If $\ba_{\bullet}$ is a graded sequence of ideals and if
$\bb_{\bullet}$ is the corresponding sequence of asymptotic multiplier ideals,
then $\lc(\ba_{\bullet})=\lc(\bb_{\bullet})$.
\end{theorem}

\begin{proof}
For every $m$, we have $\ba_m\subseteq\bb_m$, hence $\lc(\ba_m)\leq \lc(\bb_m)$.
Multiplying by $m$ and taking the limit, gives $\lc(\ba_{\bullet})
\leq \lc(\bb_{\bullet})$.

On the other hand, for fixed $m$, let $p$ be such that
$\bb_m=\cI((1/p)\cdot\ba_{mp})$. Lemma~\ref{bound_lc} below
gives
$$\frac{1}{\lc(\bb_m)}\geq\frac{1}{p\cdot \lc(\ba_{mp})}-1.$$
Dividing by $m$ and using $mp\cdot \lc(\ba_{mp})\leq \lc(\ba)$, gives
$$\frac{1}{m\cdot \lc(\bb_m)}\geq \frac{1}{\lc(\ba_{\bullet})}-\frac{1}{m}.$$
Taking the limit when $m$ goes to infinity, gives the other inequality
that we need. 
\end{proof}

\begin{lemma}\label{bound_lc}
For every non-zero ideal $\ba\subseteq R$, and every $\lambda>0$, we have
$$\frac{1}{\lc(\cI(\lambda\cdot\ba))}\geq\frac{\lambda}{\lc(\ba)}-1.$$
\end{lemma}

\begin{proof}
We prove first the following general fact: for every ideal $\ba\subseteq R$,
and every $\lambda$, $\mu>0$, we have
$$(\star)\,\,\,\,\overline{\cI(\mu\cdot\cI(\lambda\cdot\ba))
^{1/(\mu+1)}}\subseteq\cI\left(
\frac{\lambda\mu}{\mu+1}\cdot\ba)\right).$$
Recall that for an ideal $I\subseteq R$ and for $\alpha>0$, 
$$\overline{I^{\alpha}}=\{u\in R\mid{\rm ord}_E(u)\geq\alpha\cdot
{\rm ord}_E(I),
\text{for all}\,E\},$$
where $E$ ranges over all divisors over $X=\Spec\,R$.

To prove $(\star)$, let $u$ be an element in the left hand side, 
and let $E$ be a divisor over
$X$. We pick a smooth model $X'$ on which $E$ is a divisor and denote
by $K$ the relative canonical divisor of $X'$ over $X$.
By the definition
of multiplier ideals, we have
$${\rm ord}_E(u)> (1/(\mu+1))\cdot
\left(\mu\cdot\left(\lambda\cdot{\rm ord}_E(\ba)-
{\rm ord}_E(K)-1\right)-{\rm ord}
_E(K)-1\right).$$
An easy computation gives 
$${\rm ord}_E(u)>\frac{\lambda\mu}{\mu+1}\cdot{\rm ord}_E(\ba)-
{\rm ord}_E(K)-1,$$ 
hence $u\in\cI(\lambda\mu/(\mu+1)\cdot\ba)$.

It follows from $(\star)$ that if $\mu< \lc(\cI(\lambda\cdot\ba))$,
then $\lambda\mu<(\mu+1)\cdot \lc(\ba)$. Since we may assume $\lambda> \lc(\ba)$
(otherwise the statement of the lemma is trivial), we deduce
$\lc(\cI(\lambda\cdot\ba))\leq \lc(\ba)/(\lambda-\lc(\ba)),$
which immediately gives the assertion in the lemma.
\end{proof}

\begin{remark}
There are other invariants that one can associate to $\ba_{\bullet}$
and $\bb_{\bullet}$.  For example, fix a divisor $E$ over $\Spec\,R$.
Then the sequence of numbers 
$\{{\rm ord}_E(\ba_m)\}_m$ satisfies the hypothesis
in Lemma~\ref{limit}, hence we may define
$${\rm ord}_E(\ba_{\bullet}):=\lim_{m\to\infty}\frac{{\rm ord}_E(\ba_m)}{m}
=\inf_m\frac{{\rm ord}_E(\ba_m)}{m}.$$ 
Similarly, by Lemma~\ref{limit1},
we may define
$${\rm ord}_E(\bb_{\bullet}):=\lim_{m\to\infty}\frac{{\rm ord}_E(\bb_m)}{m}
=\sup_m\frac{{\rm ord}_E(\bb_m)}{m}.$$

It is easy to show that ${\rm ord}_E
(\ba_{\bullet})={\rm ord}(\bb_{\bullet})$.
Indeed, since $\ba_m\subseteq\bb_m$, we have ${\rm ord}_E(\bb_m)\leq
{\rm ord}_E(\ba_m)$, hence ${\rm ord}_E(\bb_{\bullet})
\leq {\rm ord}_E(\ba_{\bullet})$.

For the reverse inequality, fix a model $X'$ over $X=\Spec\,R$
on which $E$ is a divisor, and let $K$ be the relative canonical divisor
of $X'/X$. It follows from the definition of multiplier ideals that
$${\rm ord}_E(\bb_m)>{\rm ord}_E(\ba_m)-{\rm ord}_E(K)-1.$$
Dividing by $m$ and taking the limit gives ${\rm ord}_E(\bb_{\bullet})
\geq {\rm ord}_E(\ba_{\bullet})$.

Consider, for example, the case when $E$ is the exceptional divisor
of the blowing-up of $X$ at the maximal ideal ${\bf m}$. For any ideal
$\ba$ of $R$, we have ${\rm ord}_E(\ba)=\max\{p\mid\ba\subseteq {\bf m}^p\}$.
In this case, ${\rm ord}_E(\ba_{\bullet})$ is denoted by $\nu(\ba_{\bullet})$
and is called the Lelong number of $\ba_{\bullet}$.

Note that if $\ba_m$ is zero-dimensional for all $m$, then
we clearly have $\ee(\ba_{\bullet})\geq \nu(\ba_{\bullet})^n$ 
and Theorem~\ref{main}
implies that $\ee(\ba_{\bullet})=0$ if and only if $\nu(\ba_{\bullet})=0$. 
\end{remark}

\smallskip

The following theorem gives an inequality beween the multiplicity
and the log canonical threshold of a 
graded sequence of zero-dimensional ideals. In the case of one ideal, this 
appeared in \cite{dfelm}, generalizing the corresponding inequality
due to Corti, for the case of surfaces (see \cite{corti}).
Generalizing from one ideal 
to a graded sequence is straightforward. However,
the main point is that our 
results on graded sequences can be used to simplify
the proof even in the case of one ideal. Note that the proof in
\cite{dfelm} also used deformation to monomial ideals, but needed
a more careful analysis of the monomial case, to get a similar inequality
between the length and the log canonical threshold.

\begin{theorem}
If $\ba$ is a graded sequence of zero-dimensional ideals in $R$, then
$$\ee(\ba_{\bullet})\geq n^n/\lc(\ba_{\bullet})^n.$$
\end{theorem}

\begin{proof}
It is enough to prove that for every zero-dimensional ideal
$I\subseteq R$, we have $\ee(I)\geq n^n/\lc(I)^n$.
Indeed, if we apply this inequality for $\ba_m$, divide by $m^n$
and take the limit, we get the assertion of the theorem.

Since $R$ is smooth, it is standard to reduce the problem to an ideal
in a polynomial ring. We may therefore suppose that
$I$ is an ideal in $R=K[X_1,\ldots,X_n]$ which is
supported at the origin ($K$ might not be algebraically closed,
but this does not cause any problems). 
We  first assume that we know the inequality in the case of a monomial ideal.

Fix a monomial order and apply
Corollary~\ref{reduction_of_computation} to get $\ee(I)$ in terms 
of multiplicities of monomial ideals:
$\ee(I)=\lim_{m\to\infty}\ee(\lin_>(I^m))/m^n$. On the other hand,
it follows from the semicontinuity property of log canonical thresholds
(see \cite{demailly} or \cite{mustata}) that 
$\lc(I)/m=\lc(I^m)\geq \lc(\lin_>(I^m))$. Since we have
$$\ee(\lin_>(I^m))\geq n^n/\lc(\lin_>(I^m))\geq m^nn^n/\lc(I)^n,$$
it is enough to divide by $m^n$ and take the limit. 

We have therefore reduced the assertion
to the case when $I$ is a monomial ideal.
In this case we have the following direct argument that we have learned
from Lawrence Ein. 

Let $P=P_I$ be the Newton polyhedron of $I$ and $c=\lc(I)$. We know that
$(1/c)\cdot(1,\ldots,1)$ lies on the boundary of $P$. Fix a facet of
$P$ with equation $\sum_iX_i/a_i=1$, which contains this point. 
We therefore have $c=\sum_i1/a_i$.
On the other hand, we have $\ee(I)=n!\cdot {\rm vol}(P)\geq \prod_ia_i$.
Therefore the inequality between the arithmetic and the geometric mean
of the numbers $\{1/a_i\}_i$ gives our inequality.
\end{proof}

\providecommand{\bysame}{\leavevmode \hbox \o3em {\hrulefill} \thinspace}

\end{document}